\documentclass[10pt,a4paper]{article}
\usepackage[cp1251]{inputenc}
\usepackage{amsfonts, amssymb, amsthm, amsmath}
\usepackage{mathtext}
\usepackage{longtable}
\usepackage[russian]{babel}
\usepackage{color}

\usepackage[left=1cm,right=1cm,top=2cm,bottom=2cm]{geometry}
\usepackage{graphicx}

\theoremstyle{remark}

\renewcommand{\r}{\mathbb R}

\renewcommand{\le}{\leqslant}
\renewcommand{\ge}{\geqslant}
\newcommand{\I}{\mathbb{I}}

\renewcommand{\phi}{\varphi}

\newcommand{\eqd}{\stackrel{d}{=}}
\newcommand{\tod}{\stackrel{d}{\longrightarrow}}

\renewcommand{\refname}{References}
\title{Modeling high-frequency order flow imbalance by functional limit
theorems for two-sided risk processes\thanks{Research supported by
Russian Scientific Foundation, project 14-11-00364.}}
\author{V. Yu. Korolev\footnote{Faculty of
Computational Mathematics and Cybernetics, Lomonosov Moscow State
University; Institute of Informatics Problems, Russian Academy of
Sciences; victoryukorolev@yandex.ru}, A. V. Chertok\footnote{Faculty
of Computational Mathematics and Cybernetics, Lomonosov Moscow State
University, Euphoria Group LLC; a.v.chertok@gmail.com}, A. Yu.
Korchagin \footnote{Faculty of Computational Mathematics and
Cybernetics, Lomonosov Moscow State University;
sasha.korchagin@gmail.com} A. I. Zeifman\footnote{Vologda State
University; Institute of Informatics Problems, Russian Academy of
Sciences; a$\_$zeifman@mail.ru}}
\date{}

\begin{document}

\maketitle

{\small

{\bf Abstract:} A micro-scale model is proposed for the evolution of
the limit order book. Within this model, the flows of orders
(claims) are described by doubly stochastic Poisson processes taking
account of the stochastic character of intensities of bid and ask
orders that determine the price discovery mechanism in financial
markets. The process of {\it order flow imbalance} (OFI) is studied.
This process is a sensitive indicator of the current state of the
limit order book since time intervals between events in a limit
order book are usually so short that price changes are relatively
infrequent events. Therefore price changes provide a very coarse and
limited description of market dynamics at time micro-scales. The OFI
process tracks best bid and ask queues and change much faster than
prices. It incorporates information about build-ups and depletions
of order queues so that it can be used to interpolate market
dynamics between price changes and to track the toxicity of order
flows. The {\it two-sided risk processes} are suggested as
mathematical models of the OFI process. The multiplicative model is
proposed for the stochastic intensities making it possible to
analyze the characteristics of order flows as well as the
instantaneous proportion of the forces of buyers and sellers, that
is, the intensity imbalance (II) process, without modeling the
external information background. The proposed model gives the
opportunity to link the micro-scale (high-frequency) dynamics of the
limit order book with the macro-scale models of stock price
processes of the form of subordinated Wiener processes by means of
functional limit theorems of probability theory and hence, to give a
deeper insight in the nature of popular subordinated Wiener
processes such as generalized hyperbolic L{\'e}vy processes as
models of the evolution of characteristics of financial markets. In
the proposed models, the subordinator is determined by the evolution
of the stochastic intensity of the external information flow.

\smallskip


}

\section{Introduction}

High-frequency trading became a significant portion of the trading
volume on equity, futures and options exchanges. The high-frequency
behavior of the so-called limit order book, the list of trading
orders, is now a popular object of stochastic modeling, see, e.g.,
\cite{Parlour1998}, \cite{Foucault1999}, \cite{Goettler_et_al2005},
\cite{Avellaneda_Stoikov2008}, \cite{Rosu2009}. R. Cont et al.
\cite{ContRamaStoikov2010b} proposed a continuous-time Markov model
for the limit order book dynamics. In the paper
\cite{ContRamaStoikov2010b} the limit order book is considered as a
special queuing system where incoming orders and cancelations of
existing orders of unit sizes arrive according to independent
Poisson processes. Such kind of queuing system can be described in
terms of birth-death processes, where the states represent the
number of shares at a given order price, and transitions take place
by birth (the entry of a new limit order), or death (removal from
the limit order book by cancelation or matching with a new market
order). Birth-death processes are well-studied statistical models
which can be regarded as special examples of more general {\it
two-sided risk processes}, the stochastic models known in insurance
mathematics as {\it risk processes with stochastic premiums}.

A queueing-system-type mathematical model for the limit order book
\cite{ContRamaStoikov2010b} is supplied with rather strict formal a
priori conditions. For instance, one of such conditions is that the
intensities of order flows are assumed constant. On the one hand,
these assumptions provide the possibility to calculate at least some
characteristics of the limit order book dynamics. However, on the
other hand, these assumptions turn out to be too restrictive and
unrealistic from the practical point of view. Therefore, it is
extremely desirable to have a convenient integral characteristic of
the current state of the limit order book which can be calculated
and studied without the queueing systems framework.

Such a characteristic, the {\it order flow imbalance $($OFI$)$
process}, was introduced in 2011 in the paper \cite{Cont2011}, the
final version of which \cite{Cont2014} was published in 2014. The
same process was independently introduced and studied in
\cite{Korolev_et_al_Isr, Korolev_et_al_2013, Korolev_et_al_2013a,
Chertok2014} under the name of {\it generalized price process}. The
OFI process appears to be considerably more sensitive to the market
information than the price process itself. This is due to that time
intervals that are involved in modern high-frequency trading
applications are usually so short that price changes are relatively
infrequent events. Therefore price changes provide a very coarse and
limited description of market dynamics. However, OFI takes account
of not only best ask and bid changes, but also the
arrivals/cancelations of orders deep inside the limit order book,
which is of interest and importance because each such event
influences the current balance between buyers and sellers and
fluctuates on a much faster timescale than prices. It incorporates
information about build- ups and depletions of order queues and it
can be used to interpolate market dynamics between price changes.

In the framework of the approach developed in the present paper the
principal idea is that the moments at which the state of the limit
order book changes form a {\it chaotic} point stochastic process on
the time axis. Moreover, this point process turns out to be
non-stationary (time-non-homogeneous) because the changes of the
state of the limit order book are to a great extent subject to the
influence of non-stationary information flows. As is known, most
reasonable probabilistic models of non-stationary
(time-non-homogeneous) chaotic point processes are {\it doubly
stochastic Poisson processes} also called {\it Cox processes} (see,
e. g., \cite{Grandell1976, BeningKorolev2002}). These processes are
defined as Poisson processes with stochastic intensities. Pure
Poisson processes can be regarded as best models of stationary
(time-homogeneous) chaotic flows of events \cite{BeningKorolev2002}.
Recall that the attractiveness of a Poisson process as a model of
homogeneous discrete stochastic chaos is due to at least two
reasons. First, Poisson processes are point processes characterized
by that time intervals between successive points are independent
random variables with one and the same exponential distribution and,
as is well known, the exponential distribution possesses the maximum
differential entropy among all absolutely continuous distributions
concentrated on the nonnegative half-line with finite expectations,
whereas the entropy is a natural and convenient measure of
uncertainty. Second, the points forming the Poisson process are
uniformly distributed along the time axis in the sense that for any
finite time interval $[t_1,t_2]$, $t_1<t_2$, the conditional joint
distribution of the points of the Poisson process which fall into
the interval $[t_1,t_2]$ under the condition that the number of such
points is fixed and equals, say, $n$, coincides with the joint
distribution of the order statistics constructed from an independent
sample of size $n$ from the uniform distribution on $[t_1,t_2]$
whereas the uniform distribution possesses the maximum differential
entropy among all absolutely continuous distributions concentrated
on finite intervals and very well corresponds to the conventional
impression of an absolutely unpredictable random variable (see, e.
g., \cite{GnedenkoKorolev1996, BeningKorolev2002}).

Financial markets are examples of complex open systems whose
behavior is subject to randomness of two types: internal
(endogenous) and external (exogenous). The internal source produces
the uncertainty due to the difference of the strategies of a very
large number of traders. The <<physical>> analog of such a
randomness is the chaotic thermal motion of particles in closed
systems. The external source of randomness is a poorly predictable
flow of political or economical news influencing the interests and
strategies of traders. These two sources of randomness will be taken
into account when we construct the models of order flow imbalance
process and number of orders imbalance process as two-sided risk
processes driven by a Cox process, that is, by a Poisson process
with stochastic intensity.

As is known, the empirical (statistical) distributions of increments
of (the logarithms of) financial indices and, in particular, of
stock prices on comparatively short time intervals have
significantly heavy-tailed distributions with vertices noticeably
sharper than those of normal distributions, that is, they are {\it
leptokurtic}. At the same time, as is observed in some studies,
financial time series possess the fractal property, that is, they
are to some extent self-similar at different time scales, see, e.
g., \cite{Mandelbrot1997}.

Stable L{\'e}vy processes were among the first models applied in
practice successfully explaining the observed leptokurticity of
finite-dimensional distributions of the processes in financial
markets as well as their self-similarity. According to the approach
based on classical limit theorems of probability theory, non-normal
stable L{\'e}vy processes can appear as limits in functional limit
theorems for random walks only if the variances of elementary jumps
are infinite. In \cite{ContLarrard2011} some diffusive-type
functional limit theorems were proved for the dynamics of limit
order books in liquid markets and it was especially indicated that
within the approach used in that paper it is possible to obtain
stable limit processes only if the order sizes posses infinite
variances. In turn, the latter is possible only if the probabilities
of arbitrarily large jumps are positive. Unfortunately, the latter
condition looks very doubtful from the practical viewpoint.
Therefore, within the framework of the classical approach used in
\cite{Mandelbrot1997} and \cite{ContLarrard2011}, the theoretical
explanation of stable L{\'e}vy processes as adequate models for the
evolution of stock prices or financial indexes is at least
questionable.

In financial mathematics the evolution of (the logarithms of) stock
prices and financial indexes on small time horizons is often modeled
by random walks. The simplest example of such an approach is the
Cox--Ross--Rubinstein model (see, e. g., \cite{Shiryaev1999}). At
the same time most successful (adequate) models of the dynamics of
(the logarithms of) financial indexes on large time horizons are
subordinated Wiener processes (processes of Brownian motion with
random time) such as generalized hyperbolic processes, in
particular, variance gamma (VG) processes and
normal$\backslash\!\backslash$inverse Gaussian (NIG) processes, see
\cite{Shiryaev1999}. Subordinated Wiener processes have
finite-dimensional distributions possessing the properties mentioned
above: these distributions are heavy-tailed and leptokurtic.

In \cite{Clark1973, Korolev1994} and \cite{GnedenkoKorolev1996}
the heavy-tailedness of the empirical distributions of stock
prices was explained with the use of limit theorems for sums of a
{\it random number} of independent random variables as a
particular case of randomly stopped random walks.

When random walks are considered, the scheme of random summation is
a natural analog of the scheme of subordination of more general
random processes. In \cite{GnedenkoKorolev1996, Korolev2000,
Korolev2011} it was proposed to model the evolution of
non-homogeneous chaotic stochastic processes, in particular, of the
dynamics of stock prices and financial indexes, by random walks
generated by compound doubly stochastic Poisson processes (compound
Cox pocesses). A {\it doubly stochastic Poisson process} (also
called a {\it Cox process}) is a stochastic point process of the
form $N_1(\Lambda(t))$, where $N_1(t)$, $t\geq0$, is a homogeneous
Poisson process with unit intensity and the stochastic process
$\Lambda(t)$, $t\geq0$, is independent of $N_1(t)$ and possesses the
following properties: $\Lambda(0)=0$, ${\sf P}(\Lambda(t)<\infty)=1$
for any $t>0$, the sample paths of $\Lambda(t)$ do not decrease and
are right-continuous. A compound Cox process is a random sum of
independent identically distributed random variables in which the
number of summands follows a Cox process. Similar continuous-time
random walks were considered in \cite{gz97, z95, Huang2012}.

In accordance with the approach used in \cite{GnedenkoKorolev1996,
Korolev2000, Korolev2011} the limit distributions for compound Cox
process with elementary jumps possessing finite variances must have
the form of scale-location mixtures of normal laws which are always
heavy-tailed and leptokurtic, if the mixing distribution is
non-degenerate. Moreover, in \cite{Korolev1997, Korolev1998} it was
shown that non-normal stable laws can appear as limit distributions
for sums of independent identically distributed random variables
with {\it finite} variances, if the number of summands in the sum is
random and its distribution converges to a stable law concentrated
on the nonnegative half-line. In terms of compound Cox processes the
latter condition means that finite-dimensional distributions of the
leading process $\Lambda(t)$ are asymptotically stable. In turn,
this means that the intensity of the flow of informative events is
essentially irregular which results in the well-known clustering
effect.

In \cite{KorolevZaksZeifman2013a} some functional limit theorems
were proved establishing convergence of random walks generated by
compound Cox processes with jumps possessing finite variances to
L{\'e}vy processes with symmetric distributions including symmetric
strictly stable L{\'e}vy processes. Here we extend these results to
a non-symmetric case and apply them to modeling the evolution of the
OFI process, an integral characteristic of the behavior of the limit
order book.

Functional limit theorems are a quite natural link between random
walks and subordinated Wiener processes. The operation of
subordination gives a good explanation of the presence of heavy
tails in the empirical distributions of the increments of (the
logarithms of) stock prices and financial indexes. The functional
limit theorems for compound Cox processes proved in
\cite{KorolevZaksZeifman2013a} serve as a bridge between formal
micro-scale models having the form of continuous-time random walks
generated by compound Cox processes and popular macro-scale models
of the form of subordinated Wiener processes including generalized
hyperbolic processes, variance gamma processes, etc. The practical
importance of these models is justified by the stochastic character
of the intensities of chaotic flows of informative events in large
financial information systems and, in particular, in high-frequency
trading systems. The use of high-frequency statistical data
available due to electronic trading systems makes it possible to
verify the models mentioned above and to link them with the price
discovery process resulting from the limit order book evolution.

This paper presents a further development of the models and
techniques proposed in our previous papers \cite{Korolev_et_al_Isr},
\cite{Korolev_et_al_2013}, \cite{Chertok2014}. In this paper we
propose a convenient and rather realistic model for the OFI process
using the notion of {\it two-sided risk processes}. For this purpose
we use the techniques developed in \cite{KorolevZaksZeifman2013a}.

The paper is organized as follows. We construct our model step by
step. In Section 2 we give necessary definitions and under the
condition that the intensities of order flows are constant introduce
the two-sided risk process adopting the notion of a risk process
with stochastic premiums from insurance mathematics and prove that
it is a special compound Poisson process. In Section 3 we introduce
the conditional non-homogeneous OFI process, introduce the
multiplicative representation for the intensities of order flows in
the limit order book and, finally, introduce the assumption that the
intensity of the external information flow is stochastic and obtain
the general (unconditional) OFI process as a special compound doubly
stochastic Poisson process. The remaining part of the paper is
devoted to functional limit theorems for the so-defined OFI process.
Section 4 contains some preliminary material on the Skorokhod space
and L{\'e}vy processes. In section 5 we prove general functional
limit theorem establishing the conditions for convergence of OFI
processes to L{\'e}vy processes in the Skorokhod space in terms of
the behavior of the intensities of order flow. For this purpose we
slightly extend the classical results presented, say, in
\cite{JacodShiryaev2003}. In Section 6 we consider the conditions
for the convergence of OFI processes with elementary jumps (that is,
order sizes) possessing finite variances to the L{\'e}vy processes
with variance-mean mixed normal one-dimensional distributions, that
is, to subordinated Wiener processes, in particular, to generalized
hyperbolic L{\'e}vy processes. In Section 7 we discuss empirical
results and the final adjustment of the model.

\section{Basic model of order book dynamics and conditional
homogeneous order flow imbalance process}

Market participants trading an asset in an order-driven market can
undertake three different kinds of actions. They may
\begin{itemize}
\item[$(a)$] place a {\it limit order} to buy or sell a specified number of
shares of the asset at a particular price specified at the time of
the order,
\item[$(b)$] place a {\it market order} to buy or sell a specified number of shares
of the asset at the best currently available price, which is
executed immediately, or
\item[$(c)$] {\it cancel} a previously placed limit order that has not yet been executed.
\end{itemize}
The outstanding limit orders are summarized in a {\it limit order
book}, which lists the total number of shares of buy and sell limit
orders at each price. This limit order book is constantly changing
as new orders arrive, and we are interested in modeling statistical
properties of the limit order book state. The limit buy orders are
called {\it bids}, and the limit sell orders are called {\it asks}.
The size of the order is the number of shares specified by the
order. The lowest price for which there is an outstanding limit sell
order is called the {\it best ask} and the highest price for which
there is an outstanding limit buy order is called the {\it best
bid}. The gap between the best ask and the best bid is called the
{\it bid-ask} spread. The average of the best ask and best bid is
called the {\it mid-price}.

Of course, a limit order can turn turn out to be a market order, if
the price specified in it allows immediate matching with one of the
limit orders on the opposite side of the limit order book.


We consider the dynamics of the limit order book on the discrete
price lattice $\Pi = \{ 1, 2, ..., M \}$ as a continuous-time
process
$$
\mathbf{x}(t) \equiv (\mathbf{V^a}(t); \mathbf{V^b}(t)) \equiv
(V^a_1(t), V^a_2(t), ..., V^a_M(t); V^b_1(t), V^b_2(t), ...,
V^b_M(t)),\ \ \ t\ge 0,
$$
where $V^a_p(t)$ ($V^b_p(t))$ is the number of sell (buy) limit
orders with price $p \in \Pi$. Since buy and sell orders with one
and the same specified price cannot exist at a time (otherwise they
are immediately matched), we necessarily have $V^a_p(t) \vee
V^b_p(t) = 0$ for all $p\in\Pi$ and $t$.

The best ask $a(t)$ is defined as $a(t) = \inf \{p: V^a_p(t) >
0\}\wedge (M + 1)$, the best bid $b(t)$ is $ b(t) = \sup \{p:
V^b_p(t) > 0\} \vee 0$. Correspondingly, the mid-price is defined as
$P(t) = \frac12[a(t) + b(t)]$. Thus, the price process $P(t)$ is a
result of the evolution of the limit order book generated by the
flow of orders of three types.

First let us assume that the information flow coming from the outer
medium is fixed. Then, under fixed information, we can assume that
the inner randomness has a stable chaotic character. As it has
already been said in Sect. 1, Poisson processes are natural
mathematical models of continuous-time chaotic point processes
characterized by that the intervals between informative events
(arrivals of orders) are independent identically distributed random
variables with exponential distribution, see, e. g.,
\cite{GnedenkoKorolev1996, Korolev2011, KorolevBeningShorgin2011}.
Therefore, on the first stage of the construction of our model the
order flows are assumed to be independent renewal processes with
exponentially distributed interrenewal times (as it was done, e. g.,
in the papers {\cite{ContRamaStoikov2010b}},
{\cite{ContLarrard2011}}):

\begin{itemize}
\item limit buy (sell) orders independently arrive at a price level situated at the
distance $i$ from the best ask (bid) so that the time intervals
between successive arrivals have exponential distribution with
parameter $\lambda_i^{+}>0$ $(\lambda_i^{-}>0)$ (empirical studies
{\cite{ZovkoFarmer2002}} and {\cite{Bouchaud2002}} show that the
power law $\lambda_i^{\pm} = ki^{-\alpha}$ is a good approximation);
\item market buy (sell) orders arrive independently so that the time intervals
between successive arrivals have exponential distribution with
parameter $\mu^{+}$ $(\mu^{-})$;
\item cancelations of a limit buy (sell) order at a price level situated at the
distance $i$ from the best bid (ask) arrive with frequency
$\theta_i^{+}$ $(\theta_i^{-})$.
\end{itemize}

Under the above assumptions, $\mathbf{x}(t)$ is a continuous-time
Markov chain with state space $(\mathbb{Z}^{+})^{2M}$ and the
transitions
\begin{align*}
V^a_i(t) \to V^a_i(t) + 1 & \text{ with intensity } \lambda^{-}_{i - b(t)} & \text{ for } i > b(t), \\
V^a_i(t) \to V^a_i(t) - 1 & \text{ with intensity } \theta^{-}_{i - a(t)} & \text{ for } i \ge a(t), \\
V^a_i(t) \to V^a_i(t) - 1 & \text{ with intensity } \mu^{+} & \text{ for } i = a(t) > 0. \\
V^b_i(t) \to V^b_i(t) + 1 & \text{ with intensity } \lambda^{+}_{a(t) - i} & \text{ for } i < a(t), \\
V^b_i(t) \to V^b_i(t) - 1 & \text{ with intensity } \theta^{+}_{b(t) - i} & \text{ for } i \le b(t), \\
V^b_i(t) \to V^b_i(t) - 1 & \text{ with intensity } \mu^{-} & \text{
for } i = b(t) < M + 1.
\end{align*}
Thus, we can define the following independent Poisson processes:
\begin{itemize}
\item $L^{\pm}_{i}(t):$ the flows of limit orders with intensities $\lambda_{i}^{\pm}$;
\item $M^{\pm}(t):$ the flows of market orders with intensities
$\mu^{+}\mathbb{I}(\mathbf{V^a} \ne 0)$ and
$\mu^{-}\mathbb{I}(\mathbf{V^b} \ne 0)$;
\item $C^{\pm}_{i}(t):$ the flows of cancelations of limit orders with intencities $\theta_{i}^{\pm}$,
\end{itemize}
and the Poisson process
$$
N(t) = M^{+}(t) + M^{-}(t) + \sum\nolimits_{i = 1}^{M} (L_i^{+}(t) +
L_i^{-}(t)) + \sum\nolimits_{i = 1}^{M} (C_i^{+}(t) + C_i^{-}(t)),
$$
which describes the flow of all arriving orders.

The processes $L^{\pm}_{i}(t)$, $M^{\pm}(t)$, $C^{\pm}_{i}(t)$
completely determine the price process $P(t)$ for which the
corresponding stochastic differential equations can be written out
\cite{AbergelJedidi2011}. However, its further analytical
interpretation is very difficult or even impossible even under
rather strong and unrealistic assumption that the intensities of the
flows of orders of different types are constant.

It seems reasonable to consider some indicator of the current state
of the limit order book taking account of not only best ask and bid
changes, but also the arrivals/cancelations of orders deep inside
the limit order book, which is of interest and importance because
each such event influences the current balance of the forces of
buyers and sellers.

Recall that for the time being the intensity of the external
information flow is assumed to be fixed (constant). First fix a time
interval $[0, T]$ which is rather small so that the intensities of
the events described above can be assumed constant within this
interval. Let, as above, $N(t)$, $t\in[0,T]$, be the Poisson process
counting all the events in the limit order book and having the
intensity
$$
\lambda = \mu^{+} + \mu^{-} + \sum\nolimits_{i = 1}^{M}
(\lambda_i^{+} + \lambda_i^{-}) + \sum\nolimits_{i = 1}^{M}
(\theta_i^{+} + \theta_i^{-}).
$$
Split $N(t)$ into the sum of two independent Poisson processes
$N^+(t)$ and $N^-(t)$ with the corresponding intensities
$$
\lambda^+ = \mu^{+} + \sum\nolimits_{i=1}^M \lambda_i^{+} +
\sum\nolimits_{i=1}^M \theta_i^{-}
$$
and
$$
\lambda^- = \mu^{-} + \sum\nolimits_{i=1}^M \lambda_i^{-} +
\sum\nolimits_{i=1}^M \theta_i^{+}.
$$
Thus, $\lambda = \lambda^+ + \lambda^-$ and the processes $N^+(t)$
and $N^-(t)$ characterize the cumulative force of buyers and
sellers, correspondingly (note that a cancelation of sell orders
increases the force of buyers and vice versa). The processes
$N^+(t)$ and $N^-(t)$ are conditionally independent under fixed
outer information flow during the time interval $[0, t]$.

Following \cite{Korolev_et_al_2013a}, consider the following model
of a conditional homogeneous OFI process. Let $X_1^+, X_2^+, ...$ be
identically distributed non-negative random variables corresponding
to sequential buy order sizes and also let $X_1^-, X_2^-, ...$ be
identically distributed nonnegative random variables corresponding
to sequential sell order sizes. Let $N_1^-(t)$ and $N_1^+(t)$ be two
standard Poisson processes (that is, homogeneous Poisson processes
with unit intensities). Assume that for each $t$, the random
variables $X_1^+, X_2^+, ...$, $X_1^-, X_2^-, ...$, $N_1^-(t)$ and
$N_1^+(t)$ are independent. Introduce the process
$$
Q^{(CH)}(t)=\sum\nolimits_{j=1}^{N_1^+(\lambda^+t)}X_j^+-\sum\nolimits_{j=1}^{N_1^-(\lambda^-t)}X_j^-.\eqno(1)
$$
The process $Q^{(CH)}(t)$ will be called the {\it conditional
homogeneous OFI process}. This process is an integral instantaneous
characteristic of the state of the limit order book under the ideal
condition that the intensity of the external information flow does
not change.

Formally, the process $Q^{(CH)}(t)$ introduced above is nothing else
but what is well known in insurance mathematics as {\it the risk
process with stochastic premiums}, the positive term of which
describes the flow of insurance premiums while the negative term
describes the flow of insurance claims, see, e. g.,
\cite{KorolevBeningShorgin2011}. However, as it was especially noted
in \cite{KorolevBeningShorgin2011}, such models can hardly be
assumed adequate in insurance practice because the two components of
$Q^{(CH)}(t)$ are independent whereas in real insurance practice the
counting process in the negative component of $Q^{(CH)}(t)$ is a
rarefaction of that in the positive component and, hence, the
components cannot be assumed stochastically independent in insurance
practice. Nevertheless, as it was also noted in
\cite{KorolevBeningShorgin2011}, model (1) can be successfully
applied to the description of the processes of speculative financial
activity. In financial applications, processes of type (1) describe
the balance of forces of buyers and sellers and the resulting risks.
Hence, in what follows we will use a terminology more appropriate
for the financial applications and call the processes of form (1)
{\it two-sided risk processes}.

\smallskip

{\sc Lemma 1.} {\it Under the above assumptions, $Q^{(CH)}(t)$ is a
compound Poisson process. Namely, if $N_1(t)$ is the standard
Poisson process, then for each $t\ge0$
$$
Q^{(CH)}(t)\eqd\sum\nolimits_{j=1}^{N_1((\lambda^++\lambda^-)\,t)}X_j,\eqno(2)
$$
where $X_1,X_2,...$ are identically distributed random variables
with the common characteristic function
$$ \mathfrak f (s)\equiv{\sf E}e^{isX_1}
=\frac{\lambda^+\mathfrak
f^+(s)}{\lambda^++\lambda^-}+\frac{\lambda^-\mathfrak
f^-(-s)}{\lambda^++\lambda^-},\ \ \ s\in\mathbb{R},\eqno(3)
$$
where $\mathfrak{f}^+(s)$ and $\mathfrak{f}^-(s)$ are the
characteristic functions of the random variables $X_1^+$ and
$X_1^-$, respectively. Moreover, for each $t\ge0$ the random
variables $N_1\big((\lambda^++\lambda^-)\,t\big)$, $X_1,X_2,...$ are
independent.}

\smallskip

{\sc Proof}. Since $Q^{(CH)}(t)$ is a difference of two independent
homogeneous processes with independent increments, it also possesses
these properties. So, it remains to consider its characteristic
function. We obviously have
$$
\mathfrak f_{Q^{(CH)}(t)}(s)\equiv{\sf
E}e^{isQ^{(CH)}(t)}=\exp\big\{\lambda^+t\big(\mathfrak{f}^+(s)-1\big)\big\}\exp\big\{\lambda^-t\big(\mathfrak{f}^-(-s)-1\big)\big\}=
$$
$$
\exp\bigg\{(\lambda^++\lambda^-)t\bigg[\frac{\lambda^+\mathfrak
f^+(s)}{\lambda^++\lambda^-}+\frac{\lambda^-\mathfrak
f^-(-s)}{\lambda^++\lambda^-}\bigg]-1\bigg\},\ \ \ s\in\mathbb{R}.
$$
The latter characteristic function corresponds to the compound
Poisson distribution specified in the formulation of the lemma.

\smallskip

{\sc Remark 1}. It is easy to see that the random variable $X_1$ is
a randomization:
$$
X_1=\begin{cases}X_1^+&\text{with probability }
{\displaystyle\frac{\lambda^+}{\lambda^++\lambda^-}},\vspace{2mm}\\
-X_1^-&\text{with probability }
{\displaystyle\frac{\lambda^-}{\lambda^++\lambda^-}}\end{cases}
$$
so that
$$
{\sf E}X_1=\frac{\lambda^+{\sf
E}X_1^+}{\lambda^++\lambda^-}-\frac{\lambda^-{\sf
E}X_1^-}{\lambda^++\lambda^-},
$$
$$
{\sf D}X_1=\bigg(\frac{\lambda^+}{\lambda^++\lambda^-}\bigg)^2{\sf
D}X_1^++\bigg(\frac{\lambda^-}{\lambda^++\lambda^-}\bigg)^2{\sf
D}X_1^--\frac{\lambda^+\lambda^-}{(\lambda^++\lambda^-)^2}{\sf
E}\big(X_1^+-X_1^-\big)^2.
$$
Some empirical studies (see, e. g., \cite{Chakraborti_et_al2009,
Toke_2010}) show that new limit order sizes are in fact randomly
distributed according to the exponential law resulting in that
within model (2), in general, the distribution of independent random
variables $X_j$ is the asymmetric Laplace law (a mixture of the
exponential distribution and the distribution symmetric to
exponential concentrated on the negative half-line).

\section{General order flow imbalance process}

In this section the ideal conditional model described above will be
adapted with the account of non-constant character of the intensity
of external information flow resulting in that the parameters
$\lambda^+$ and $\lambda^-$ describing the corresponding reaction of
buyers and sellers vary in time.

In reality, the intensities of the order flows are non-homogeneous
since the external information flow determining the intensities of
the events in the order book is non-homogeneous itself. In this case
the intensities of the flows of orders of different types, first,
may not be independent and, second, depend in a certain way on some
process which determines the external news background. To formalize
these ideas, assume that the intensities of the processes introduced
in Sect. 2 vary in time:
$$
\mu^+=\mu^+(t),\ \ \lambda^+_i=\lambda^+_i(t),\ \ \
\theta^+_i=\theta^+_i(t),
$$
$$
\mu^-=\mu^-(t),\ \ \lambda^-_i=\lambda^-_i(t),\ \ \
\theta^-_i=\theta^-_i(t),
$$
$i=1,\ldots,M$. For $t\ge0$ introduce positive functions
$$
\lambda^+(t)=\mu^+(t)+\sum\nolimits_{i=1}^M\lambda^+_i(t)+\sum\nolimits_{i=1}^M\theta^+_i(t),
$$
$$
\lambda^-(t)=\mu^-(t)+\sum\nolimits_{i=1}^M\lambda^-_i(t)+\sum\nolimits_{i=1}^M\theta^-_i(t),
$$
and introduce the functions
$$
\Lambda^+(t)=\int_0^t\lambda^+(\tau)d\tau\ \ \text{and}\ \
\Lambda^-(t)=\int_0^t\lambda^-(\tau)d\tau, \ \ t\ge0.
$$
Let $N^+_1(t)$ and $N^-_1(t)$ be two independent Poisson processes,
each having the unit intensity. Put
$$
N^+(t)=N^+_1\big(\Lambda^+(t)\big),\ \ \
N^-(t)=N^-_1\big(\Lambda^-(t)\big).
$$
The processes $N^+_1(t)$ and $N^-_1(t)$ are non-homogeneous
Poisson processes with the instantaneous intensities
$\lambda^+(t)$ and $\lambda^-(t)$, respectively. The process
$$
Q^{(CN)}(t)=\sum\nolimits_{j=1}^{N^+(t)}X_j^+-\sum\nolimits_{j=1}^{N^-(t)}X_j^-
$$
will be called the {\it conditional non-homogeneous generalized
price process}. This process is an integral instantaneous
characteristic of the state of the limit order book under the
condition that the intensity of the external information flow varies
non-randomly.

Empirical data analyzed in \cite{Korolev_et_al_2013a} give very
serious grounds to assume that actually both $\Lambda^+(t)$ and
$\Lambda^-(t)$ depend on {\it one and the same} process
$\Lambda^*(t)$ describing the general agitation of the market as its
reaction on the flow of external information, see Figure
\ref{fig:lambda_est} below. So, in the subsequent reasoning we will
assume that
$$
\Lambda^+(t)=\alpha^+(t)\Lambda^*(t),\ \ \ \
\Lambda^-(t)=\alpha^-(t)\Lambda^*(t)
$$
where $\alpha^+(t)$ and $\alpha^-(t)$ are some positive functions.
Moreover, it was shown in \cite{Korolev_et_al_2013a} that the above
multiplicative representations for the cumulative intensities
$\Lambda^+(t)$ and $\Lambda^-(t)$ with the application of special
limit theorems directly lead to the conclusion that the asymptotic
(<<heavy-traffic>>) approximations for the statistical regularities
of the behavior of OFI must have the form of normal variance-mean
mixtures. An extremely high adequacy of normal variance-mean
mixtures, in particular, of generalized hyperbolic laws, as models
of statistical regularities of the behavior of characteristics of
financial markets noted in, say, the canonical works
\cite{CarrMadanChang1998, EberleinKeller1995,
EberleinKellerPrause1998, EberleinPrause1998, Eberlein1999,
MadanSeneta1990, Prause1997, Shiryaev1999} can serve as a
theoretical evidence in favor of the multiplicativity representation
for the intensities of flows of informative events in financial
markets. In Section 6 we will discuss this topic in more detail.

For simplicity we further assume that $\alpha^+(t)$ and
$\alpha^-(t)$ are constant so that
$$
\alpha^+(t)\equiv\alpha^+>0,\ \ \ \alpha^-(t)\equiv\alpha^->0.
$$
According to \cite{Korolev_et_al_2013a}, the process $\Lambda^*(t)$
can be interpreted as the amplifying factor of the trading
intensities due to a very uncertain and poorly predictable external
news. Therefore, concerning the function $\Lambda^*(t)$, we make a
somewhat more general assumption. In what follows we will suppose
that $\Lambda^*(t)$ is a random measure, that is, it is a stochastic
process possessing the following properties: $\Lambda^*(0)=0$, ${\sf
P}(\Lambda^*(t)<\infty)=1$ for any $t>0$, the sample paths of
$\Lambda^*(t)$ do not decrease and are right-continuous. Moreover,
assume that the process $\Lambda^*(t)$ is independent of the
standard Poisson processes $N^+(t)$ and $N^-(t)$. In what follows
the process
$$
Q(t)=\sum\nolimits_{j=1}^{N^+(\alpha^+\Lambda^*(t))}X^+_j-\sum\nolimits_{j=1}^{N^-(\alpha^-\Lambda^*(t))}X^-_j
$$
will be called the (general) {\it order flow imbalance $($OFI$)$
process}.

Under the conditions just imposed, in accordance with Lemma 1 we
have
$$
{\sf P}\big(Q(t)<x\big)=\int_0^{\infty}{\sf
P}\bigg(\sum\nolimits_{j=1}^{N^+(\alpha^+\lambda)}X^+_j-\sum\nolimits_{j=1}^{N^-(\alpha^-\lambda)}X^-_j<x\bigg)d{\sf
P}\big(\Lambda^*(t)<\lambda\big)=
$$
$$
=\int_0^{\infty}{\sf
P}\bigg(\sum\nolimits_{j=1}^{N_1((\alpha^++\alpha^-)\lambda)}X_j<x\bigg)d{\sf
P}\big(\Lambda^*(t)<\lambda\big)={\sf
P}\bigg(\sum\nolimits_{j=1}^{N_1((\alpha^++\alpha^-)\Lambda^*(t))}X_j<x\bigg),\eqno(4)
$$
where $X_1,X_2,...$ are identically distributed random variables
with the common characteristic function
$$
\mathfrak f (s)\equiv{\sf E}e^{isX_1} =\frac{\alpha^+\mathfrak
f^+(s)}{\alpha^++\alpha^-}+\frac{\alpha^-\mathfrak
f^-(-s)}{\alpha^++\alpha^-},\ \ \ s\in\mathbb{R},\eqno(5)
$$
and $N_1(t)$ is the standard Poisson process, moreover, all the
processes and random variables involved in representation (4) are
independent.

Denote $\Lambda(t)=(\alpha^++\alpha^-)\Lambda^*(t)$. Obviously, the
process $\Lambda(t)$ is a random measure independent of the standard
Poisson process $N_1(t)$ taking part in representation (4). Let
$N(t)=N_1\big(\Lambda(t)\big)$. $N(t)$ is a doubly stochastic
Poisson process (Cox process). Based on (4), in what follows by an
OFI process we will mean the process
$$
Q(t)=\sum\nolimits_{j=1}^{N_1(\Lambda(t))}X_j,
$$
where the random variables $X_1,X_2,...$ have the common
characteristic function (5) and all the involved random variables
and processes are assumed independent. The process $Q(t)$ so defined
is a special two-sided risk process whose positive and negative
components are not independent: they are linked by one and the same
process describing the cumulative intensity of the flows of
<<positive>> and <<negative>> events in the limit order book.

This process is a sensitive indicator of the current state of the
limit order book since time intervals between events in a limit
order book are usually so short that price changes are relatively
infrequent events. Therefore price changes provide a very coarse and
limited description of market dynamics at time micro-scales. The OFI
process tracks best bid and ask queues and change much faster than
prices. It incorporates information about build-ups and depletions
of order queues so that it can be used to interpolate market
dynamics between price changes and to track the toxicity of order
flows.

\begin{figure}[!h]
  \centerline{
    \includegraphics[width=0.5 \textwidth, height=0.35\textwidth]{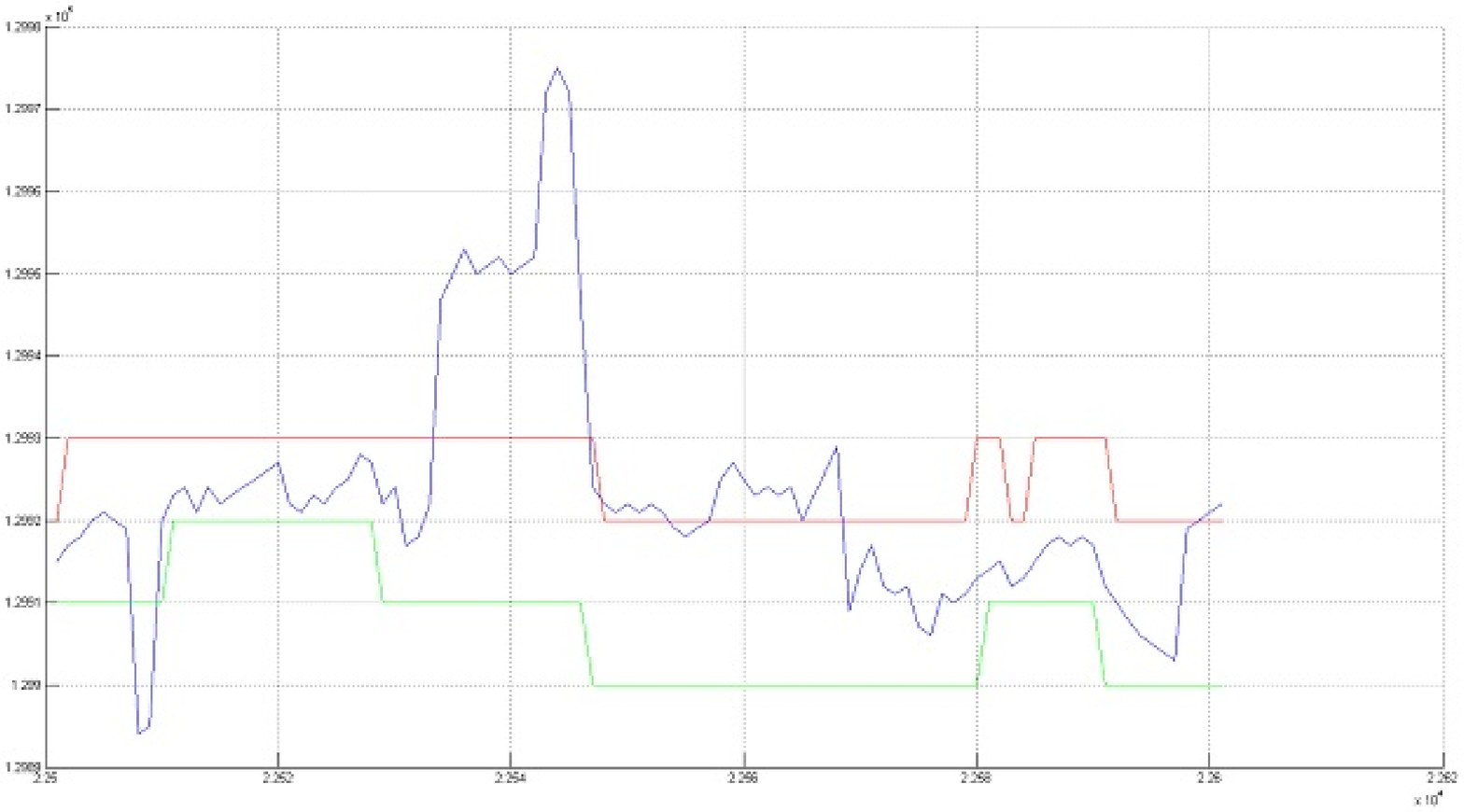}
    \includegraphics[width=0.5 \textwidth,
    height=0.35\textwidth]{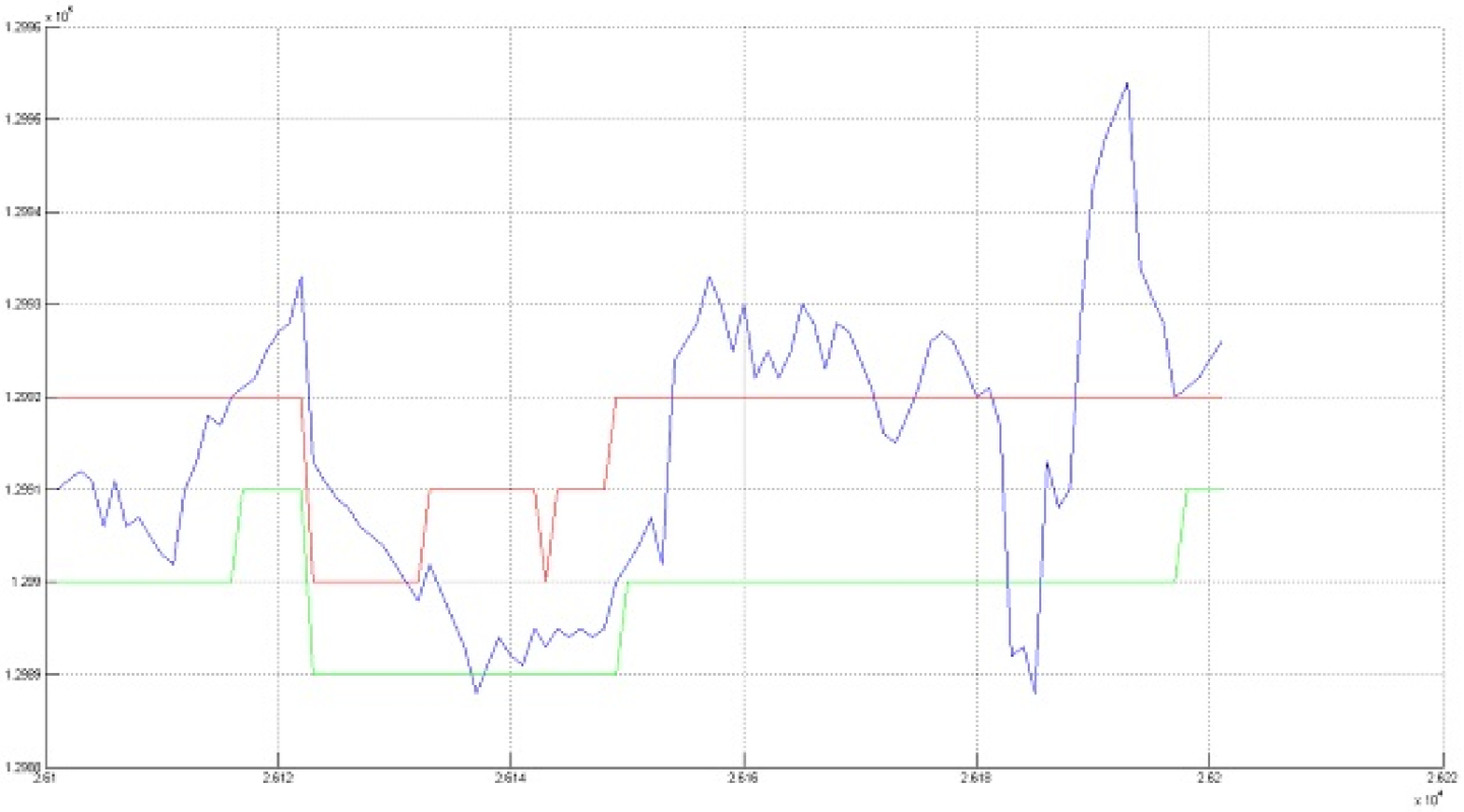}}
  \caption{The price (horizontal lines) and the OFI within two
  100-milliseconds time intervals on 01/07/2014, futures of RTSI.}\label{fig:adv_sel}
\end{figure}

Figure \ref{fig:adv_sel} well illustrates that the OFI process is
noticeably more sensitive and informative than the price process
since it is obviously more volatile. On this figure the horizontal
lines mark the values of the price for the RTS index and broken line
corresponds to the OFI process within two 100-milliseconds time
intervals.

The OFI process $Q(t)$ in some sense accumulates more information
than the pure price process which can be constructed in the same way
so that the flows of <<positive>> and <<negative>> events contain
only those events each of which results in the change of the current
price upward and downward respectively. In this sense the price
process is a kind of rarefaction of the OFI process and can be
studied by the same techniques.

\smallskip

The above assumptions give the opportunity to use the well-developed
analytic apparatus of compound Cox processes to study the asymptotic
behavior of the OFI process under the condition that the intensities
of the flows of informative events related to the limit order book
are large. This will enable us to describe possible asymptotic
approximations to the OFI process. To do so, we need some auxiliary
definitions and results.

\section{Skorokhod space. L{\'e}vy processes}

Let $D=D[0,1]$ be the space of real functions defined on $[0,1]$,
right-continuous and having finite left-side limits.

Let $\mathcal{F}$ be the class of strictly increasing mappings of
$[0,1]$ onto itself. Let $f$ be a non-decreasing function on $[0,1]$
with $f(0)=0$, $f(1)=1$. Set $\|f\|=\sup_{s\neq
t}\left|\log\left[\big(f(t)-f(s)\big)/(t-s)\right]\right|$. If
$\|f\|<\infty$, then the function $f$ is continuous and strictly
increasing and, hence, belongs to $\mathcal{F}$.

Define the distance $d_0(x,y)$ in the set $D[0,1]$ as the greatest
lower bound of the set of positive numbers $\epsilon$, for which
$\mathcal{F}$ contains a function $f$ such that $\|f\|\le\epsilon$
and $\sup_t|x(t)-y(f(t))|\le\epsilon$.

It can be shown that the space $D[0,1]$ is complete with respect to
the distance $d_0$. The metric space $\mathcal{D}=(D[0,1],d_0)$ is
called {\it the Skorokhod space}. Everywhere in what follows we will
consider stochastic processes as $\mathcal{D}$-valued random
elements.

Let $X, X_1,X_2,...$ be $\mathcal{D}$-valued random elements. Let
$T_X$ be a subset of $[0,1]$ such that $0\in T_X$, $1\in T_X$ and if
$0<t<1$, then $t\in T_X$ if and only if ${\sf P}\bigl(X(t)\neq
X(t-)\bigr)=0$. The following theorem establishing sufficient
conditions for the weak convergence of stochastic processes in
$\mathcal{D}$ (denoted below as $\Longrightarrow$ and assumed as
$n\to\infty$) is well-known.

\smallskip

{\sc Theorem A.} {\it Let}
$\bigl(X_n(t_1),...,X_n(t_k)\bigr)\Longrightarrow
\bigl(X(t_1),...,X(t_k)\bigr)$ {\it for any natural $k$ and
$t_1,...,t_k$ belonging to $T_X$. Let ${\sf P}\bigl(X(1)\neq
X(1-)\bigr)=0$ and let there exist a non-decreasing continuous
function $F$ on $[0,1]$, such that for any $\epsilon>0$}
$$
{\sf P}\bigl(|X_n(t)-X_n(t_1)|\ge\epsilon,\
|X_n(t_2)-X_n(t)|\ge\epsilon\bigr)
\le\epsilon^{-2\nu}\bigl[F(t_2)-F(t_1)\bigr]^{2\gamma}\eqno(6)
$$
{\it for $t_1\le t\le t_2$ and $n\ge1$, where $\nu\ge 0$,
$\gamma>1/2$. Then} $X_n\Longrightarrow X$.

\smallskip

The {\sc proof} of Theorem A can be found, for example, in
\cite{Billingsley1968}.

\smallskip

Everywhere in what follows the symbol $\eqd$ stands for the
coincidence of distributions.

By a L{\'e}vy process we will mean a stochastic process $X(t)$,
$t\ge0$, possessing the following properties: (i) $X(0)=0$ almost
surely; (ii) $X(t)$ is the process with independent increments, that
is, for any $N\ge1$ and $t_0,t_1,...,t_N$ ($0\le t_0\le t_1\le...\le
t_N$) the random variables $X(t_0)$, $X(t_1)-X(t_0)$, $...,$
$X(t_N)-X(t_{N-1})$ are jointly independent; (iii) $X(t)$ is a
homogeneous process, that is, $X(t+h)-X(t)\eqd X(s+h)-X(s)$ for any
$s,t,h>0$; (iv) the process $X(t)$ is stochastically continuous,
that is, for any $t\ge0$ and $\epsilon>0$ $ \lim_{s\to t}{\sf
P}(|X(t)-X(s)|>\epsilon)=0$; (v) sample paths of the process $X(t)$
are right-continuous and have finite left-side limits.

Denote the characteristic function of the random variable $X(t)$ as
$\psi_t(s)$ ($\psi_t(s)={\sf E}e^{isX(t)}$, $s\in\mathbb{R}$). The
following statement describes a well-known property of L{\'e}vy
processes.

\smallskip

{\sc Lemma 2.} {\it Let $X=X(t)$, $t\ge 0$, be a L{\'e}vy process.
For any $t>0$ the characteristic function of the random variable
$X(t)$ is infinitely divisible and has the form
$$
\psi_{t}(s) = \left[\psi_{1}(s)\right]^{t} = \left[{\sf
E}\,e^{isX(1)}\right]^t,\ \ \ \  s\in\mathbb{R}.\eqno(7)
$$
Conversely, let $Y$ be an arbitrary infinitely divisible random
variable. Then the family of infinitely divisible distributions with
characteristic functions of the form $\left[{\sf
E}\,e^{isY}\right]^t$ completely determines finite-dimensional
distributions of a L{\'e}vy process $X(t)$, $t\ge 0$, moreover,}
$X(1)\stackrel{d}{=}Y$.

\smallskip

The properties of L{\'e}vy processes are described in detail in
\cite{Bertoin1996, Sato1999}. The books
\cite{BarndorffNielsenMikoschResnick2001, Schoutens2003} and the
review paper \cite{Geman2002} deal with applications of L{\'e}vy
processes to modeling the dynamics of stock prices and financial
indexes.

By $G_{\alpha,\theta}(x)$ we will denote the distribution function
of the strictly stable law with the characteristic exponent $\alpha$
and parameter $\theta$ corresponding to the characteristic function
$\mathfrak{g}_{\alpha,\theta}(s)=\exp\big\{-|s|^{\alpha}\exp\big\{-\frac{i}{2}
\pi\theta\alpha\mathrm{sign}s\big\}\big\}$, $s\in\r$, where
$0<\alpha\le2$,
$|\theta|\le\theta_{\alpha}=\min\{1,\frac{2}{\alpha}-1\}$. To
symmetric strictly stable distributions there corresponds the value
$\theta=0$. To one-sided stable distributions there correspond the
values $\theta=1$ and $0<\alpha\le1$.

If $\xi$ is a random variable with the distribution function
$G_{\alpha,\theta}(x)$, $0<\alpha<2$, then ${\sf
E}|\xi|^{\delta}<\infty$ for any $\delta\in(0,\alpha)$, but the
moments of orders higher or equal to $\alpha$ of the random variable
$\xi$ do not exist (see, e. g., \cite{Zolotarev1983}).

The distribution function of the standard normal law ($\alpha=2$,
$\theta=0$) will be denoted $\Phi(x)$,
$\Phi(x)=\int_{-\infty}^x\phi(z)dz,$
$\phi(x)=\frac{1}{\sqrt{2\pi}}e^{-x^2/2}$.

It is well known that the distribution function $G_{\alpha,0}(x)$ of
the symmetric strictly stable law can be represented as a scale
mixture of normal laws:
$$
G_{\alpha,0}(x)=\int_{0}^{\infty}\Phi\Big(\frac{x}{\sqrt{u}}\Big)dG_{\alpha/2,1}(u),\
\ \ x\in\mathbb{R}\eqno(8)
$$
(see, e.g., \cite{Zolotarev1983}, Theorem 3.3.1). To representation
(8) there corresponds the analogous representation in terms of
characteristic functions:
$$
\mathfrak{g}_{\alpha,0}(s)=\int_{0}^{\infty}\exp\Big\{-\frac{s^2u}{2}\Big\}dG_{\alpha/2,1}(u),\
\ \ s\in\mathbb{R}.
$$

A L{\'e}vy process $X(t)$, $t\ge0$, will be called {\it
$\alpha$-stable}, if ${\sf P}\big(X(1)<x\big)=G_{\alpha,\theta}(x)$,
$x\in\mathbb{R}$. It can be shown (see, e.g.,
\cite{EmbrechtsMaejima2002}) that if $X(t)$, $t\ge0$, is a L{\'e}vy
process, then $X(t)$ is $\alpha$-stable if and only if
$$
X(t)\eqd t^{1/\alpha}X(1),\ \ \ t\ge0.
$$

\section{Convergence of OFI processes to L{\'e}vy processes}

In what follows without noticeable loss of generality we will
consider stochastic processes defined for $0\le t\le1$. Actually,
this means that we consider the behavior of generalized two-sided
risk processes and hence, OFI processes, on finite time horizons.
The equality of the right bound of the horizon to one can be
achieved by an appropriate choice of the units of measurement of
time. In other words, we will concentrate on studying the case of
the Skorokhod space $\mathcal{D}$.

In order to introduce reasonable asymptotics which formalizes the
condition of <<infinite>> growth of intensities of the flows of
informative events, and makes it possible to construct asymptotic
(<<heavy-traffic>>) approximations to the one-dimensional
distributions of the OFI process, fix a time instant $t$ and
introduce an auxiliary parameter $n$. Everywhere in what follows the
convergence will be meant as $n\to\infty$ unless otherwise
specified. So, consider a sequence of compound Cox processes of the
form
$$
Q_n(t)=\sum\nolimits_{i=1}^{N^{(n)}_{1}(\Lambda_n(t))}X_{n,i},\ \ \
t\ge0,\eqno(9)
$$
where $\{N^{(n)}_{1}(t),\, t\ge0\}_{n\ge1}$ is a sequence of Poisson
processes with unit intensities; for each $n=1,2,...$ the random
variables $X_{n,1},X_{n,2},...$ are identically distributed; for any
$n\ge1$ the random variables $X_{n,1},X_{n,2},...$ and the process
$N^{(n)}_{1}(t)$, $t\ge0$, are independent; for each $n=1,2,...$
$\Lambda_n(t)$, $t\ge0$, is a subordinator, that is, a
non-decreasing positive L{\'e}vy process, independent of the process
$$
Z_n(t)=\sum\nolimits_{i=1}^{N^{(n)}_{1}(t)}X_{n,i},\ \ \ t\ge
0,\eqno(10)
$$
and such that $\Lambda_n(0)=0$ and there exist $\delta\in(0,1]$,
$\delta_1\in(0,1]$ and the constants $C_n\in(0,\infty)$ providing
for all $t\in(0,1]$ the validity of the inequality
$$
{\sf E}\Lambda^{\delta}_n(t)\le (C_nt)^{\delta_1}.\eqno(11)
$$
Here and in what follows for definiteness we assume
$\sum\nolimits_{i=1}^0=0$. In terms introduced in section 2,
$Z_n(t)$ is a conditional homogeneous OFI process.

Below we will demonstrate that the reasoning presented in
\cite{KorolevZaksZeifman2013a} for the symmetric case can be by
slight changes generalized to the case of limit processes with
non-symmetric distributions.

From (9) and (10) it is easy to see that $Q_n(t)=Z_n(\Lambda_n(t))$.
Since for each $n\ge1$ both $Z_n(t)$ and $\Lambda_n(t)$ are
independent L{\'e}vy processes, and, moreover, $\Lambda_n(t)$ is a
subordinator, then the superposition $Q_n(t)=Z_n(\Lambda_n(t))$ is
also a L{\'e}vy process (see, e. g., Theorem 3.1.1 in
\cite{Kashcheev2001}). Hence the following statement follows.

\smallskip

{\sc Lemma 3}. {\it For any $0\le t_1<t_2<\infty$ and any $n\ge1$ we
have} $Q_n(t_2)-Q_n(t_1)\eqd Q_n(t_2-t_1)$.

\smallskip

Denote $a_n={\sf E}X_{n,1}$ and assume that
$$
0<m_n^{\beta}\equiv{\sf E}|X_{n,1}|^{\beta}<\infty\eqno(12)
$$
for some $\beta\in[1,2]$.

\smallskip

{\sc Remark 3}. If we supply the parameters introduced in section 3
by the index $n$ and look at the topic under consideration here from
the viewpoint of modeling the dynamics of the limit order book, then
we can assume that
$$
{\sf E}e^{isX_{n,1}} =\frac{\alpha^+_n\mathfrak
f^+_n(s)}{\alpha^+_n+\alpha^-_n}+\frac{\alpha^-_n\mathfrak
f^-_n(-s)}{\alpha^+_n+\alpha^-_n},\ \ \ s\in\mathbb{R},
$$
where $\mathfrak f^+_n(s)$ and $\mathfrak f^+_n(s)$ are the
characteristic functions of the random variables $X^+_{n,1}$ and
$X^-_{n,1}$, the elementary positive and negative increments of the
OFI process, then
$$
a_n\equiv{\sf E}X_{n,1}=\frac{\alpha^+_n{\sf
E}X^+_{n,1}}{\alpha^+_n+\alpha^-_n}-\frac{\alpha^-_n{\sf
E}X^-_{n,1}}{\alpha^+_n+\alpha^-_n}
$$
and
$$
m_n^{\beta}\equiv{\sf E}|X_{n,1}|^{\beta}=\frac{\alpha^+_n{\sf
E}|X^+_{n,1}|^{\beta}}{\alpha^+_n+\alpha^-_n}+\frac{\alpha^-_n{\sf
E}|X^-_{n,1}|^{\beta}}{\alpha^+_n+\alpha^-_n}
$$
so that condition (12) is implied by the conditions
$$
0<{\sf E}|X^+_{n,1}|^{\beta}<\infty,\ \ \ 0<{\sf
E}|X^-_{n,1}|^{\beta}<\infty.
$$

\smallskip

{\sc Lemma 4}. {\it Let $Q_n(t)$ be a compound Cox process $(9)$
satisfying conditions $(11)$ and $(12)$. Then for any $t\in[0,1]$
and any $\epsilon>0$ we have} ${\sf
P}\big(|Q_n(t)|\ge\epsilon\big)\le
\epsilon^{-\beta\delta}m_n^{\beta\delta}\cdot\big(C_nt\big)^{\delta_1}$.

\smallskip

{\sc Proof}. Since one-dimensional distributions of the Cox process
(9) are mixed Poisson, we have
$$
{\sf P}\big(|Q_n(t)|\ge\epsilon\big)= 
\sum\nolimits_{k=0}^{\infty}{\sf
P}\big(N^{(n)}_1(\Lambda_n(t))=k\big){\sf
P}\Big(\Big|\sum\nolimits_{j=1}^kX_{n,j}\Big|\ge\epsilon\Big)=
$$
$$
=\int_{0}^{\infty}\Big[\sum\nolimits_{k=0}^{\infty}e^{-\lambda}\frac{\lambda^k}{k!}{\sf
P}\Big(\Big|\sum\nolimits_{j=1}^kX_{n,j}\Big|\ge\epsilon\Big)\Big]\,d{\sf
P}\big(\Lambda_n(t)<\lambda\big).\eqno(13)
$$
The change of the order of summation and integration is possible due
to the obvious uniform convergence of the series. Continue (13) by
the sequential application of the Markov and Jensen inequalities
with $\delta\in(0,1]$ taking part in (11) and $\beta\in[1,2]$ taking
part in (12). As a result we obtain
$$
{\sf P}\big(|Q_n(t)|\ge\epsilon\big)\le
\frac{1}{\epsilon^{\beta\delta}}\int_{0}^{\infty}\Big[\sum\nolimits_{k=0}^{\infty}e^{-\lambda}\frac{\lambda^k}{k!}\Big({\sf
E}\Big|\sum\nolimits_{j=1}^kX_{n,j}\Big|^{\beta}{\Big)\!}^{\delta}\,\Big]\,d{\sf
P}\big(\Lambda_n(t)<\lambda\big),\eqno(14)
$$
since with $\delta\in(0,1]$ the function $f(x)=x^{\delta}$ is
concave for $x\ge0$. It is easy to see that ${\sf
E}\big|\sum\nolimits_{j=1}^kX_{n,j}\big|^{\beta}\le\sum\nolimits_{j=1}^k{\sf
E}|X_{n,j}|^{\beta}=km_n^{\beta}$ for $1\le\beta\le2$.
Therefore, continuing (14) with the account of the Jensen inequality
for concave functions and (11), we obtain
$$
{\sf P}\big(|Q_n(t)|\ge\epsilon\big)\le
\frac{m_n^{\beta\delta}}{\epsilon^{\beta\delta}}\int_{0}^{\infty}{\sf
E}\big[N^{(1)}_1(\lambda)\big]^{\delta}\,d{\sf
P}\big(\Lambda_n(t)<\lambda\big)\le
$$
$$
\le
\frac{m_n^{\beta\delta}}{\epsilon^{\beta\delta}}\int_{0}^{\infty}\big[{\sf
E}N^{(1)}_1(\lambda)\big]^{\delta}\,d{\sf
P}\big(\Lambda_n(t)<\lambda\big) =
\frac{m_n^{\beta\delta}}{\epsilon^{\beta\delta}}\cdot{\sf
E}\Lambda_n^{\delta}(t)\le
\frac{m_n^{\beta\delta}}{\epsilon^{\beta\delta}}\cdot\big(C_nt\big)^{\delta_1}.
$$
The lemma is proved.

\smallskip

To establish weak convergence of the stochastic processes $Q_n(t)$
in the Skorokhod space $\mathcal{D}$, first it is required to find
the limit distribution of the random variables $Q_n(t)$ for each
$t>0$. The symbol $\tod$ will denote convergence in distribution,
that is, pointwise convergence of the distribution functions in all
continuity points of the limit distribution function.

Let $t=1$. Denote $N_n=N^{(n)}_{1}(\Lambda_n(1))$. Assume that for
some $k_n\in\mathbb{N}$ the convergence
$$
{\sf P}(X_{n,1}+...+X_{n,k_n}<x)\tod H(x)\eqno(15)
$$
takes place, where $H(x)$ is some infinitely divisible distribution
function.

Also assume that
$$
{\sf P}\big(\Lambda_n(1)<k_n x\big)\tod {\sf P}(U<x),\eqno(16) 
$$
where $U$ is a nonnegative random variable such that its
distribution is not degenerate in zero. Notice that since
$\Lambda_n(t)$ is a L{\'e}vy process, then the random variable $U$
is infinitely divisible being the weak limit of infinitely divisible
random variables.

\smallskip

{\sc Lemma 5}. {\it Let $N_n=N_1^{(n)}(\Lambda_n)$, $n\ge1$, where
$\{N_1^{(n)(t)},\ t\ge0\}$, $n=1,2,\ldots$ are standard Poisson
processes and $\Lambda_n$, $n=1,2,\ldots$ are positive random
variables such that for each $n\ge1$ the random variable $\Lambda_n$
is independent of the process $N_1^{(n)}(t)$. Then
$$
{\sf P}(N_n<k_nx)\tod A(x)
$$
for some infinitely increasing sequence $k_n$ of real numbers and
some distribution function $A(x)$ if and only if}
$$
{\sf P}(\Lambda_n<k_nx)\tod A(x).
$$

\smallskip

For the proof see \cite{GnedenkoKorolev1996}.

\smallskip

From Lemma 5 it follows that convergence (16) is equivalent to
$$
{\sf P}\big(N_n<k_n x\big)\tod {\sf P}(U<x).\eqno(17) 
$$
By the Gnedenko--Fahim transfer theorem \cite{GnedenkoFahim1969}
conditions (15) and (17) imply that
$$
Q_n(1)=X_{n,1}+...+X_{n,N_n}\tod Q,\eqno(18) 
$$
where $Q$ is a random variable with the characteristic function
$$
\mathfrak{f}(s)=\int_{0}^{\infty}\big(h(s)\big)^u\,d{\sf
P}(U<u),\eqno(19) 
$$
$h(s)$ being the characteristic function corresponding to the
distribution function $H(x)$. Note that the distribution function
$H(x)$  may not satisfy the condition $H(-x)=1-H(x)$ for all
$x\ge0$, that is, it may not be symmetric.

Let $Y$ be an infinitely divisible random variable with the
distribution function $H(x)$. Since both $Y$ and $U$ are infinitely
divisible, we can define independent L{\'e}vy processes $Y(t)$ and
$U(t)$, $t\ge0$, such that $Y(1)\eqd Y$ and $U(1)\eqd Y$. Then with
the account of Lemma 2 it is easy to verify that
$\mathfrak{f}(s)={\sf E}e^{isQ}={\sf E}\exp\big\{isY(U(1))\big\}$,
$s\in\mathbb{R}$, that is, $Q\eqd Y(U(1))$. Moreover, repeating the
reasoning from \cite{Kashcheev2001} (see Theorem 3.3.1 there), we
can easily see that the random variable $Q$ is infinitely divisible
and hence, we can define a L{\'e}vy process $Q(t)$, $t\ge0$, such
that $Q(1)\eqd Q$. From Lemma 2 and the abovesaid it follows that we
can regard $Q(t)$ as the superposition: $Q(t)\eqd Y(U(t))$.

Since according to (18) we have
$Q_n(1)=\sum\nolimits_{i=1}^{N_n}X_{n,i}\Longrightarrow Q(1)$, and
both $Q_n(t)$ and $Q(t)$ are L{\'e}vy processes, then, using (7) we
can conclude that for any $t>0$
$$
Q_n(t)=\sum\nolimits_{i=1}^{N_{n,1}(\Lambda_n(t))}X_{n,i}\tod Q(t).
\ \ \ \eqno(20)
$$

Since the processes $Q_n(t)$ and $Q(t)$, $0\le t\le 1$, are L{\'e}vy
processes, then almost all their sample paths belong to the
Skorokhod space ${\cal D}$.

Consider the question what additional conditions are required to
provide the weak convergence of the compound Cox process $Q_n(t)$ to
the L{\'e}vy process $Q(t)$ in the space ${\cal D}$. We will
consider each of the conditions of Theorem A one by one.

First, without loss of generality, let $0\le t_1<t_2<...<t_k \le 1$.
The convergence $\bigl(Q_n(t_1),...,Q_n(t_k)\bigr)\tod
\bigl(Q(t_1),...,Q(t_k)\bigr)$ is equivalent to the convergence
$$\bigl(Q_n(t_1),Q_n(t_2)-Q_n(t_1),...,Q_n(t_k)-Q_n(t_{k-1})\bigr)
\tod
$$
$$
\tod\bigl(Q(t_1),Q(t_2)-Q(t_1),...,Q(t_k)-Q(t_{k-1})\bigr),\eqno(21) 
$$
since the linear transform $(x_1,x_2,...,x_{k-1},x_k)\longmapsto
(x_1,x_2-x_1,...,x_k-x_{k-1})$ of $\mathbb{R}^k$ to $\mathbb{R}^k$
is one-to-one and continuous in both directions. But convergence
(21) follows from (20) and the fact that both $Q_n(t)$ and $Q(t)$
are L{\'e}vy processes.

Second, we have to check the condition ${\sf P}\bigl(Q(1)\neq
Q(1-)\bigr)=0$. This condition holds if and only if
$\lim_{t\to1-}{\sf P}\bigl(|Q(1)-Q(t)|>\epsilon\bigr)=0$ for any
$\epsilon>0$ (see relation (15.16) in \cite{Billingsley1968}).
Consider ${\sf P}\bigl(|Q(1)-Q(t)|>\epsilon\bigr)$. Since $Q(t)$ is
a L{\'e}vy process, then $Q(1)-Q(t)\eqd Q(1-t)$ by Lemma 3.
Therefore, ${\sf P}\bigl(|Q(1)-Q(t)|>\epsilon\bigr)={\sf
P}\bigl(|Q(1-t)|>\epsilon\bigr)$. For each $\epsilon>0$ and each
$t\in[0,1]$ there exists an $\epsilon_t\in[\epsilon/2,\epsilon]$
such that the points $\pm\epsilon_t$ are continuity points of the
distribution function of the random variable $Q(1-t)$. Since
$Q_n(t)\tod Q(t)$ for each $t\in[0,1]$, then ${\sf
P}\bigl(|Q(1-t)|>\epsilon_t\bigr)=\lim_{n\to\infty}{\sf
P}\bigl(|Q_n(1-t)|>\epsilon_t\bigr)$. Thus, for any $\epsilon>0$ and
any $t\in[0,1]$ we have
$$
{\sf P}\bigl(|Q(1-t)|>\epsilon\bigr)\le{\sf
P}\bigl(|Q(1-t)|>\epsilon_t\bigr)=\lim_{n\to\infty}{\sf
P}\bigl(|Q_n(1-t)|>\epsilon_t\bigr).\eqno(22) 
$$
Continuing (22) with the account of (11) and applying Lemma 4, for
$\delta\in(0,1]$ taking part in (11) we obtain
$$
{\sf P}\bigl(|Q(1-t)|>\epsilon\bigr)\le\sup_n{\sf
P}\bigl(|Q_n(1-t)|>\epsilon_t\bigr)\le
$$
$$
\le
\sup_n\big(\epsilon_t^{-\beta}m_n^{\beta}\big)^{\delta}\big(C_n|1-t|\big)^{\delta_1}\le
\big(2^{\beta\delta}\epsilon^{-\beta\delta}|1-t|\big)^{\delta_1}\sup_nm_n^{\beta\delta}C_n^{\delta_1}.\eqno(23) 
$$
Therefore, if
$$
K\equiv\sup_nC_n^{\delta_1/\delta}m_n^{\beta}<\infty,\eqno(24) 
$$
then (23) implies $\lim_{t\to1-}{\sf
P}\bigl(|Q(1)-Q(t)|\!>\!\epsilon\bigr)\!\le\!
4(K\epsilon^{-\beta})^{\delta}\lim_{t\to1-}|1-t|^{\delta_1}\!=\!0$.

Third, check condition (6) under the assumption that (11) and (24)
hold. As it has been noted above, $Q_n(t)$ is a L{\'e}vy process and
hence, it has independent increments. Therefore,
$$
{\sf
P}\big(|Q_n(t)-Q_n(t_1)|\ge\epsilon,\,|Q_n(t_2)-Q_n(t)|\ge\epsilon\big)=
$$
$$
= {\sf P}\big(|Q_n(t)-Q_n(t_1)|\ge\epsilon\big)\cdot{\sf
P}\big(|Q_n(t_2)-Q_n(t)|\ge\epsilon\big).\eqno(25) 
$$
Consider the first multiplier on the right-hand side of (25). By
Lemma 3, $Q_n(t)-Q_n(t_1)\eqd Q_n(t-t_1)$. With the account of (24),
by Lemma 4 we obtain
$$
{\sf P}\big(|Q_n(t)-Q_n(t_1)|\ge\epsilon\big)={\sf
P}\big(|Q_n(t-t_1)|\ge\epsilon\big)\le(K\epsilon^{-\beta})^{\delta}|t-t_1|^{\delta_1}.\eqno(26) 
$$
For the second multiplier on the right-hand side of (25) we
similarly obtain
$$
{\sf P}\big(|Q_n(t_2)-Q_n(t)|\ge\epsilon\big)={\sf
P}\big(|Q_n(t_2-t)|\ge\epsilon\big)\le(K\epsilon^{-\beta})^{\delta}|t_2-t|^{\delta_1}.\eqno(27) 
$$
Thus, from (26) and (27) it follows that
$$
{\sf
P}\big(|Q_n(t)-Q_n(t_1)|\ge\epsilon,\,|Q_n(t_2)-Q_n(t)|\ge\epsilon\big)\le
(K\epsilon^{-\beta})^{2\delta}\big[(t-t_1)(t_2-t)\big]^{\delta_1}\eqno(28) 
$$
It is easy to see that for any $t_1\le t\le t_2$ we have
$(t-t_1)(t_2-t)\le{\textstyle\frac14}(t_2-t_1)^2$. Substituting this
estimate in (28) we obtain ${\sf
P}\big(|Q_n(t)-Q_n(t_1)|\ge\epsilon,\,|Q_n(t_2)-Q_n(t)|\ge\epsilon\big)\le
\epsilon^{-2\beta\delta}\big[\frac12{K(t_2-t_1)}\big]^{2\delta_1}$.
Therefore, if conditions (11) and (24) hold, then condition (6)
holds with $F(t)\equiv \frac12Kt$, $\nu=\beta\delta$ and
$\gamma=\delta_1$.

Summarizing this reasoning related to checking the conditions of
Theorem A, we arrive at the following statement.

\smallskip

{\sc Theorem 1.} {\it Let the OFI processes $Q_n(t)$ $($see $(9))$
be lead by non-decreasing positive L{\'e}vy processes $\Lambda_n(t)$
satisfying conditions $(11)$ and $(16)$ with some
$\delta,\delta_1\in(0,1]$ and $k_n\in\mathbb{N}$. Assume that the
random variables $\{X_{n,j}\}_{j\ge1}$, $n=1,2,...$, $($the
randomized order sizes, that is, the jumps of the OFI process
$Q_n(t)$$)$, satisfy conditions $(15)$ with the same $k_n$ and
$(12)$ with some $\beta\in[1,2]$. Also assume that condition $(24)$
holds. Then the OFI processes $Q_n(t)$ weakly converge in the
Skorokhod space $\mathcal{D}$ to the L{\'e}vy process $Q(t)$ such
that
$$
{\sf E}\exp\{isQ(1)\}=\int_{0}^{\infty}\big(h(s)\big)^u\,d{\sf
P}(U<u),\ \ \ s\in\mathbb{R},\eqno(29) 
$$
where $h(s)$ is the characteristic function corresponding to the
distribution function $H(x)$ in $(15)$.}

\smallskip

It is worth noting that actually Theorem 1 deals with the
well-studied weak convergence of special semimartingales with
stationary increments, see, e. g., \cite{JacodShiryaev2003}.
However, the superposition-type structure of the processes
considered in the present paper makes it possible to relax the
conditions required in the general case, say, in Corollary VII.3.6
of \cite{JacodShiryaev2003} where it is assumed that (in our
terminology) $\delta=\delta_1=1$.

Some corollaries of this result dealing with symmetric limit laws
were considered in \cite{KorolevZaksZeifman2013a}. In particular, it
was demonstrated there that symmetric stable L{\'e}vy processes can
appear as limits for compound Cox processes even when the variances
of elementary increments of a compound Cox process are finite. As it
has been already said, in most applied problems there are no reasons
to reject this assumption. This is exactly the case of modeling
limit order book dynamics, where elementary increments of the OFI
process (the randomized order sizes) can be assumed bounded.
Therefore in what follows we will concentrate attention on the case
of finite variances and consider the conditions of convergence of
OFI processes to some popular models, in particular, to generalized
hyperbolic L{\'e}vy processes.

\section{Generalized hyperbolic L{\'e}vy processes as asymptotic
approximations to OFI processes}

Denote $\sigma_n^2={\sf D}X_{n,1}$. From the classical theory of
limit theorems it is well known that if, as $n\to\infty$, the
conditions
$$
k_na_n\longrightarrow a,\ \ k_n\sigma_n^2\longrightarrow \sigma^2 \
\mbox{ and } \ k_n{\sf
E}(X_{n,1}-a_n)^2\I(|X_{n,1}-a_n|\ge\epsilon)\longrightarrow
0\eqno(30) 
$$
hold for some $a\in\mathbb{R}$, $0<\sigma^2<\infty$ and any
$\epsilon>0$, then convergence (15) takes place with
$H(x)\equiv\Phi\big(\sigma^{-1}(x-a)\big)$. In this case the
distribution function $F(x)$ of the limit random variable $Q(1)$ in
Theorem 1 is a variance-mean mixture of normal laws. Recently it was
demonstrated that normal variance-mean mixtures 
appear as limiting in simple limit theorems for random sums of
independent identically distributed random variables
\cite{KorolevSokolov2012, KorolevZaks2013, Korolev2013}. Namely, let
$\{\xi_{n,j}\}_{j\ge1},$ $n=1,2,\ldots,$ be a double array of
row-wise (for each fixed $n$) independent and identically
distributed random variables. Let $\{\nu_n\}_{n\ge1}$ be a sequence
of integer nonnegative random variables such that for each $n\ge1$
the random variables $\nu_n,\xi_{n,1},\xi_{n,2},\ldots$ are
independent. Denote $ S_{n,k}=\xi_{n,1}+\ldots +\xi_{n,k}$. The
following theorem was proved in \cite{Korolev2013}.

\smallskip

{\sc Theorem B}. {\it Assume that there exist$:$ a sequence
$\{k_n\}_{n\ge1}$ of natural numbers and finite numbers
$\alpha\in\mathbb{R}$ and $\sigma>0$ such that
$$
{\sf P}\big(S_{n,k_n}<x\big)\tod
\Phi\Big(\frac{x-\alpha}{\sigma}\Big).\eqno(31)
$$
Assume that $\nu_n\to\infty$ in probability. Then the distribution
functions of random sums $S_{n,\nu_n}$ converge to some distribution
function $F(x)$$:$
$$
{\sf P}\big(S_{n,\nu_n}<x\big)\tod F(x),
$$
if and only if there exists a distribution function $A(x)$ such that
$A(0)=0$,
$$
F(x)=\int_{0}^{\infty}\Phi\Big(\frac{x-\alpha
z}{\sigma\sqrt{z}}\Big)dA(z),
$$
and
$$
{\sf P}(\nu_n<xk_n)\tod A(x).
$$
}

\smallskip

Theorem B and Lemma 5 yield the following result.

\smallskip

{\sc Theorem 2.} {\it Let the the OFI processes $Q_n(t)$ $($see
$(9))$ be lead by non-decreasing positive L{\'e}vy processes
$\Lambda_n(t)$ satisfying condition $(11)$ with some
$\delta,\delta_1\in(0,1]$. Assume that the random variables
$\{X_{n,j}\}_{j\ge1}$, $n=1,2,...$, $($the randomized order sizes,
that is, the jumps of the OFI process $Q_n(t)$$)$ satisfy conditions
$(30)$ with some $k_n\in\mathbb{N}$. Also assume that condition
$(24)$ holds with $\beta=2$. Then the OFI processes $Q_n(t)$ weakly
converge in the Skorokhod space $\mathcal{D}$ to a L{\'e}vy process
$Q(t)$ if and only if there exists a nonnegative random variable $U$
such that
$$
{\sf
P}\big(Q(1)<x\big)=\int_{0}^{\infty}\Phi\Big(\frac{x-au}{\sigma\sqrt{u}}\Big)d{\sf
P}(U<u),\ \ \ x\in\mathbb{R},\eqno(32) 
$$
and condition $(16)$ holds with the same $k_n$.}

\smallskip

The class of distributions of form (32) was systematically
considered by O. Barndorff-Nielsen and his colleagues \cite{BN1978,
BN1979, BN1982} in order to introduce {\it generalized hyperbolic
distributions} and study their properties.

The class of normal variance-mean mixtures (32) is very wide. For
example, it contains generalized hyperbolic laws with generalized
inverse Gaussian mixing distributions, in particular, $(a)$
symmetric and non-symmetric (skew) Student distributions (including
Cauchy distribution), to which in (32) there correspond inverse
gamma mixing distributions; $(b)$ variance gamma (VG) distributions)
(including symmetric and non-symmetric Laplace distributions), to
which in (32) there correspond gamma mixing distributions; $(c)$
normal$\backslash\!\backslash$inverse Gaussian (NIG) distributions
to which in (32) there correspond inverse Gaussian mixing
distributions, and many other types. Along with generalized
hyperbolic laws, the class of normal variance-mean mixtures contains
symmetric strictly stable laws with $\mu=0$ and strictly stable
mixing distributions concentrated on the positive half-line,
generalized exponential power distributions and many other types.

Generalized hyperbolic distributions demonstrate exceptionally high
adequacy when they are used to describe statistical regularities in
the behavior of characteristics of various complex open systems, in
particular, turbulent systems and financial markets. There are
dozens of dozens of publications dealing with models based on
generalized hyperbolic distributions. Just mention the canonic
papers \cite{BN1978, BN1979, BarndorffNielsen1998,
BNBlaesildSchmiegel2004, CarrMadanChang1998, EberleinKeller1995,
EberleinKellerPrause1998, EberleinPrause1998, Eberlein1999,
MadanSeneta1990, Prause1997, Shiryaev1999}. Therefore below we will
concentrate our attention on functional limit theorems establishing
the convergence of OFI processes to generalized hyperbolic L{\'e}vy
processes.

It is a convention to explain such a good adequacy of generalized
hyperbolic models by that they possess many parameters to be
suitably adjusted. But actually, it would be considerably more
reasonable to explain this phenomenon by functional limit theorems
yielding the possibility of the use of generalized hyperbolic
L{\'e}vy processes as convenient <<heavy-traffic>> {\it asymptotic}
approximations.

Denote the density of the {\it generalized inverse Gaussian
distribution} by $p_{GIG}(x;\nu,\mu,\lambda)$,
$$
p_{GIG}(x;\nu,\mu,\lambda)=\frac{\lambda^{\nu/2}}{2\mu^{\nu/2}K_{\nu}\big(\sqrt{\mu\lambda}\big)}\cdot
x^{\nu-1}\cdot\exp\Big\{-\frac12\Big(\frac{\mu}{x}+\lambda
x\Big)\Big\},\ \ \ x>0.
$$
Here $\nu\in\r$,
$$
\begin{array}{lll}
\mu>0, & \lambda\ge0, & \text{if }\nu<0,\vspace{1mm}\cr \mu>0, &
\lambda>0, & \text{if }\nu=0,\vspace{1mm}\cr \mu\ge0, & \lambda>0, &
\text{if }\nu>0,
\end{array}
$$
$K_{\nu}(z)$ is the modified Bessel function of the third kind
with index $\nu$,
$$
K_{\nu}(z)=\frac12\int_{0}^{\infty}y^{\nu-1}\exp\Big\{-\frac{z}{2}\Big(y+\frac1y\Big)\Big\}dy,\
\ \ \ z\in\mathbb{C},\ \mathrm{Re}\,z>0.
$$
The corresponding distribution function will be denoted
$P_{GIG}(x;\nu,\mu,\lambda)$,
$$
P_{GIG}(x;\nu,\mu,\lambda)=\int_{0}^{x}p_{GIG}(z;\nu,\mu,\lambda)dz,\
\ \ x\ge0,
$$
and $P_{GIG}(x;\nu,\mu,\lambda)=0$, $x<0$. According to
\cite{Seshadri1997}, the generalized inverse Gaussian distribution
was introduced in 1946 by {\'E}tienne Halphen, who used it to
describe monthly volumes of water passing through hydroelectric
power stations. In the paper \cite{Seshadri1997} generalized inverse
Gaussian distribution was called the {\it Halphen distribution}. In
1973 this distribution was re-discovered by Herbert Sichel
\cite{Sichel1973}, who used it as the mixing law in special mixed
Poisson distributions (the {\it Sichel distributions}, see, e. g.,
\cite{KorolevBeningShorgin2011}) as discrete distributions with
heavy tails. In 1977 these distributions were once more
re-discovered by O. Barndorff-Nielsen \cite{BN1977, BN1978}, who, in
particular, used them to describe the particle size distribution.

The class of generalized inverse Gaussian distributions is rather
rich and contains, in particular, both distributions with
exponentially decreasing tails (gamma-distribution ($\mu=0$,
$\nu>0$)), and distributions whose tails demonstrate power-type
behavior (inverse gamma-distribution ($\lambda=0$, $\nu<0$), inverse
Gaussian distribution ($\nu=-\frac12$) and its limit case as
$\lambda\to0$, the L{\'e}vy distribution (stable distribution with
the characteristic exponent equal to $\frac12$ and concentrated on
the nonnegative half-line, the distribution of the time for the
standard Wiener process to hit the unit level)).

In 1977--78 O. Barndorff-Nielsen \cite{BN1977, BN1978} introduced
the class of {\it generalized hyperbolic distributions} as the class
of special normal variance-mean mixtures. For convenience, we will
use a somewhat simpler parameterization. Let $\alpha\in\r$,
$\sigma>0$. If the generalized hyperbolic distribution function with
parameters $\alpha$, $\sigma$, $\nu$, $\mu$, $\lambda$ is denoted
$P_{GH}(x;\alpha,\sigma,\nu,\mu,\lambda)$, then by definition,
$$
P_{GH}(x;\alpha,\sigma,\nu,\mu,\lambda)=\int_{0}^{\infty}\Phi\Big(\frac{x-\alpha
z}{\sigma\sqrt{z}}\Big)\,p_{GIG}(z;\nu,\mu,\lambda)dz,\ \ \
x\in\r.\eqno(33)
$$
Note that in (33) mixing is carried out simultaneously by both
location and scale parameters, but since these parameters are
directly linked in (33), then actually (33) is a one-parameter
mixture. Several parameterizations are used for generalized
hyperbolic distributions, see, e. g., \cite{Shiryaev1999, BN1978,
BN1979, EberleinKeller1995, Prause1997, EberleinKellerPrause1998,
BarndorffNielsen1998, BNBlaesildSchmiegel2004}. However, the density
 $p_{GH}(x;\alpha,\sigma,\nu,\mu,\lambda)$ of the generalized
hyperbolic distribution cannot be expressed via elementary functions
and has the form
$$
p_{GH}(x;\alpha,\sigma,\nu,\mu,\lambda)=\int_{0}^{\infty}\frac{1}{\sigma\sqrt{z}}\phi\Big(\frac{x-\alpha
z}{\sigma\sqrt{z}}\Big)p_{GIG}(z;\nu,\mu,\lambda)\,dz=
$$
$$
=\int_{0}^{\infty}\frac{1}{\sigma\sqrt{2\pi
z}}\exp\Big\{-\frac{(x-\alpha
z)^2}{2\sigma^2z}\Big\}\frac{\lambda^{\nu/2}z^{\nu-1}}{2\mu^{\nu/2}K_{\nu}\big(\sqrt{\mu\lambda}\big)}\exp\Big\{-\frac12\Big(\frac{\mu
}{z}+\lambda  z\Big)\Big\}dz,\ \ \ x\in\mathbb{R},
$$
which can be further simplified using modified Bessel functions of
the third kind.

From Theorem B we easily obtain the following

\smallskip

{\sc Corollary 1}. {\it Assume that there exist$:$ a sequence
$\{k_n\}_{n\ge1}$ of natural numbers and finite numbers
$\alpha\in\mathbb{R}$ and $\sigma>0$ such that convergence $(31)$
takes place. Assume that $\nu_n\to\infty$ in probability. Then the
distribution of the random sum $S_{\nu_n}$ converges to a
generalized hyperbolic distribution$:$ ${\sf
P}\big(S_{n,\nu_n}<x\big)\tod
P_{GH}(x;\alpha,\sigma,\nu,\mu,\lambda)$, if and only if ${\sf
P}(\nu_n<xk_n)\tod P_{GIG}(x;\nu,\mu,\lambda)$. }

\smallskip

From Theorem 2 and Corollary 1 with the account of the equivalence
of relations (16) and (17) we easily obtain the following result on
the convergence of OFI processes represented as two-sided risk
processes to generalized hyperbolic L{\'e}vy processes.

\smallskip

{\sc Theorem 3.} {\it Let the OFI processes $Q_n(t)$ $($see $(9))$
be lead by non-decreasing positive L{\'e}vy processes $\Lambda_n(t)$
satisfying condition $(11)$ with some $\delta,\delta_1\in(0,1]$.
Assume that the random variables $\{X_{n,j}\}_{j\ge1}$, $n=1,2,...$,
$($the randomized order sizes, that is, the jumps of the generalized
price process $Q_n(t)$$)$ satisfy conditions $(30)$ with some
$k_n\in\mathbb{N}$ and some $a\in\mathbb{R}$ and $\sigma>0$. Also
assume that condition $(24)$ holds with $\beta=2$. Then the OFI
processes $Q_n(t)$ weakly converge in the Skorokhod space
$\mathcal{D}$ to a generalized hyperbolic L{\'e}vy process $Q(t)$
such that ${\sf
P}\big(Q(1)<x\big)=P_{GH}(x;a,\sigma,\nu,\mu,\lambda)$ if and only
if $ {\sf P}(\Lambda_n(1)<k_nx)\tod P_{GIG}(x;\nu,\mu,\lambda)$ with
the same $k_n$, $\nu$, $\mu$ and $\lambda$. }

To conclude this section we should note that Theorems 1--3 presented
above can serve as convenient explanation of the high adequacy of
generalized hyperbolic L{\'e}vy processes as models of the evolution
of OFI process. Moreover, they directly link the subordinator in the
representation of generalized hyperbolic L{\'e}vy processes as
subordinated Wiener processes with the intensities of order flows
determined by the process of general agitation of the market caused
by external news. Moreover, since the latter is hardly predictable,
the description of type of its distribution requires rather many
parameters (at least three within the family of generalized inverse
Gaussian laws).

\section{Final adjustment of the model}

Here we revisit the definition of a general order flow imbalance
process given in Section 3. In addition to the notation introduced
there, introduce the {\it intensity imbalance $($II$)$ process} as
$$
r(t)=\frac{\alpha^+(t)}{\alpha^-(t)}.
$$
Let $\Lambda^{\star}(t)=\Lambda^-(t)=\alpha^-(t)\Lambda^*(t)$. Then,
obviously, $\Lambda^+(t)=r(t)\Lambda^{\star}(t)$,
$$
\Lambda(t)=(1+r(t))\Lambda^{\star}\eqno(34)
$$
and, if $r(t)=r=\mathrm{const}$, then
$$
{\sf E}X_{n,1}=\frac{r}{1+r}{\sf E}X_{n,1}^++\frac{1}{1+r}{\sf
E}X_{n,1}^-, \ \ \ {\sf E}|X_{n,1}|^{\beta}=\frac{r}{1+r}{\sf
E}|X_{n,1}^+|^{\beta}+\frac{1}{1+r}{\sf
E}|X_{n,1}^-|^{\beta}.\eqno(35)
$$
This means that instead of using parameters $\alpha^+$ and
$\alpha^-$ we can repeat all the reasoning in terms of the intensity
imbalance process $r(t)$.

To test some of assumptions and concepts introduced above, we chose
the high-frequency data concerning trading the RTS index (RTSI)
futures, one of the most liquid financial instruments at Moscow
Exchange. We analyzed the data concerning the flows of all orders
(limit, market orders and cancelations) at the first $5$ levels of
the limit order book for the RTSI within the period from 1st to 30th
of July, 2014. These data provide an access to the most detailed
information concerning market trade unlike data concerning matchings
and quotes (TAQ, Trades and Quotes), that are often used for the
analysis of high frequency data and contain prices and volumes of
contracts (which corresponds only to market orders in the flow of
all orders), as well as the information concerning the prices and
volumes of the best bid and ask quotes (i. e., only the first level
of the limit order book) with the corresponding times. So, the stock
exchange provides full information about order flows which makes it
possible to analyze the processes $N^+(t)$ and $N^-(t)$ within the
framework of the OFI process model described above.

\renewcommand{\figurename}{{\small \bf Fig.}}

\begin{figure}[!h]
  \centerline{
    \includegraphics[width= 15cm, height=7cm]{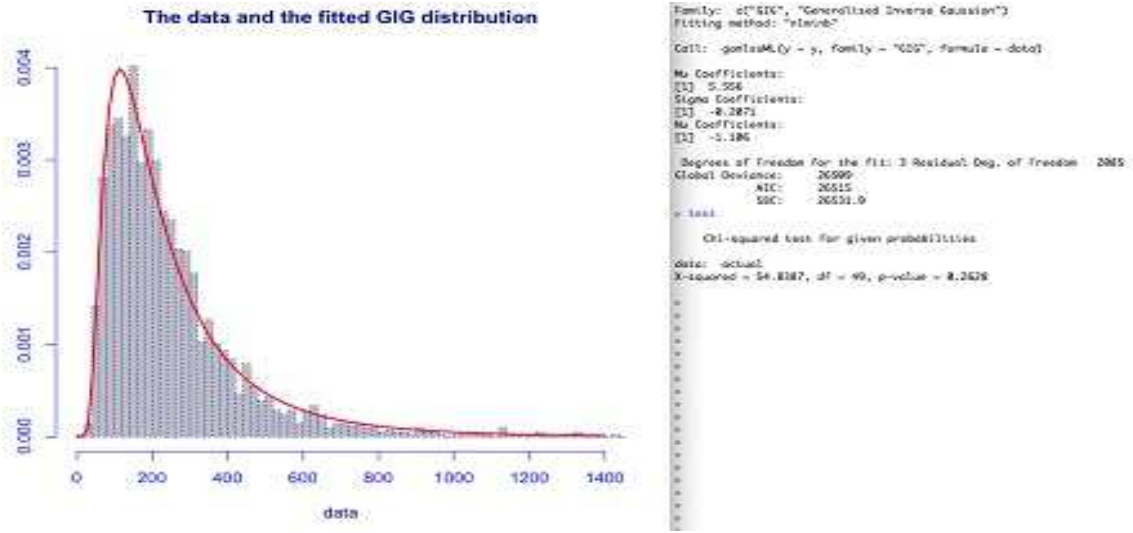}}
  \caption{\small \bf The histogram of the number of buy order arrivals per a
15-second interval. And the fitted GIG density. RTSI futures, day
session 2014.07.01).}\label{fig:nevents_hist}
\end{figure}
\begin{figure}[!h]
  \centerline{
    \includegraphics[width= 20cm, height=6cm]{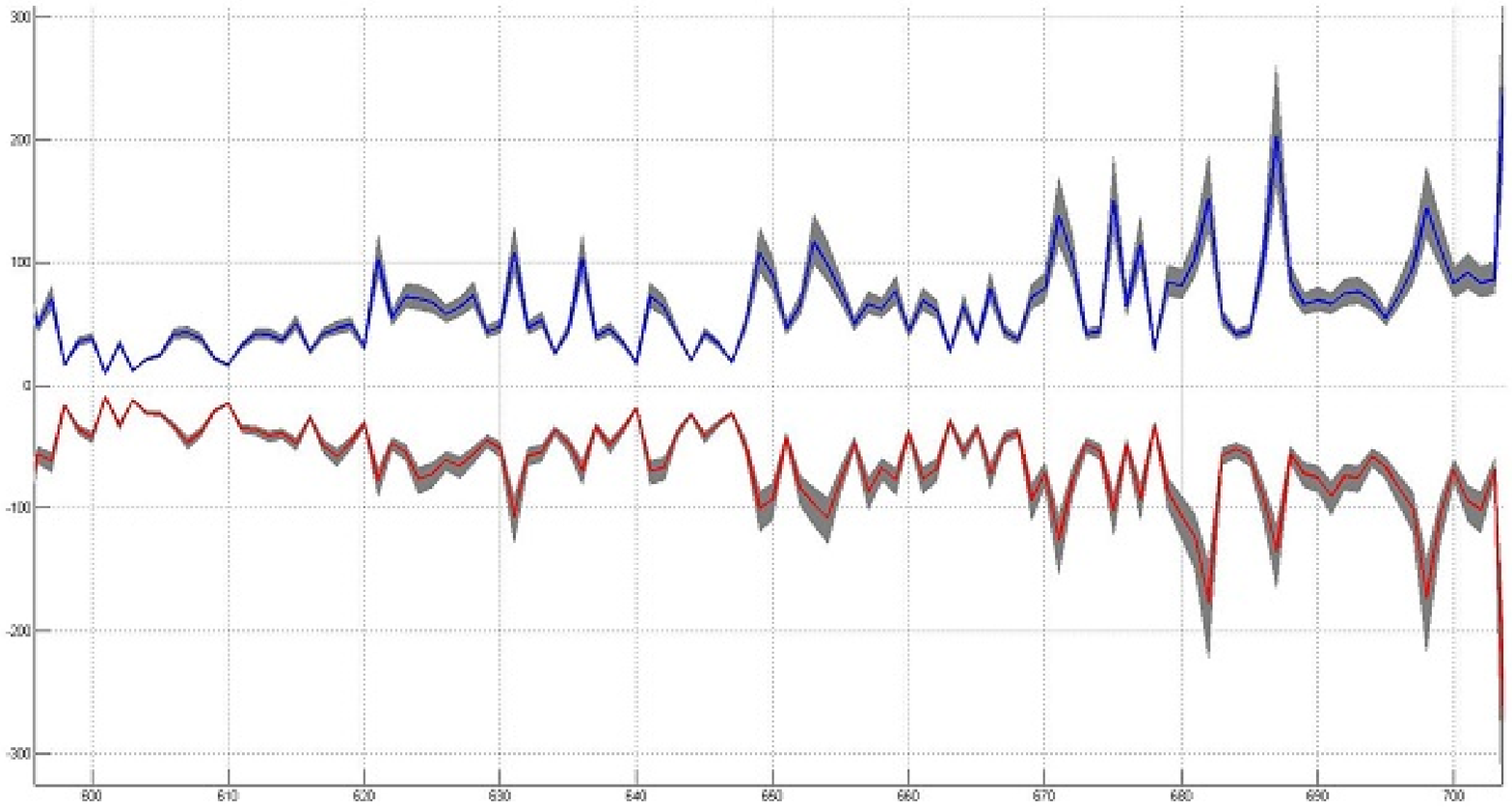}}
  \caption{\small \bf The estimates of the intensities of the buy orders flow
  (upper curve) and the mirror reflection of the same for sell
  orders flow (lower curve).}\label{fig:lambda_est}
\end{figure}
\begin{figure}[!h]
\centerline{
\includegraphics[width = 16cm, height=5cm]{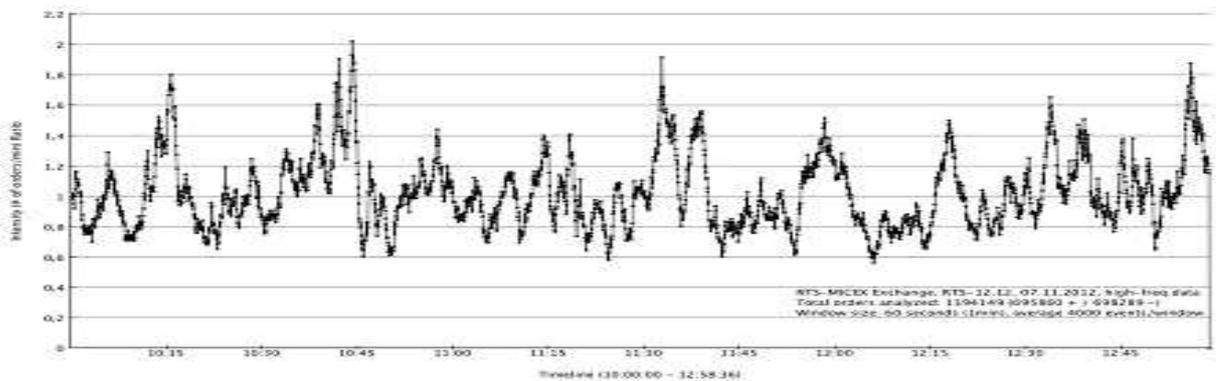}
} \caption{\small \bf  The ratio of instantaneous empirical
intensities of buyers and sellers $r(t)$ averaged over 60 sec}
\label{fig:intratio}
\end{figure}

Split one of trading days (July 1st, 2014) into time intervals of
equal 15 seconds length. We do not take into account the intervals
falling into the starting five minutes of trade (from 10:00 to
10:05), as well as those falling into the last five minutes of trade
(from 18:40 to 18:45), since these periods are characterized by
abnormal rises of the volatility hardly subject to analysis within
the framework of the proposed model.

Figure \ref{fig:nevents_hist} depicts the histogram of the number of
buy order arrivals per a 15-second interval (left) and the same for
sell order arrivals (right). It can be easily seen that these
histograms are rather satisfactorily approximated by the GIG
densities with the corresponding parameters. This is a good evidence
in favor of the Cox process model proposed in this paper for the
order counting processes. The explanation of such a good fit is
given by lemma 5: if the expected intensity is high, then the
asymptotic distribution of the mixed Poisson distribution is the
same as that of the accumulated intensity. The corresponding
intensities of the flows of buy and sell orders are shown on Fig.
\ref{fig:lambda_est} where we again see a picture convincing that
the multiplicativity representation for the intensities really
holds. On Fig. \ref{fig:intratio} the graph of estimated (averaged
over 60 sec) intensities imbalance process $r(t)$ is shown. It
illustrates the general preponderance of buyers over sellers during
the period of observation. Its deviations from the unit level are
rather small.

So, as we have seen above, if $r(t)=r=\mathrm{const}$, then the OFI
process can be successfully modeled by a generalized hyperbolic
L{\'e}vy process with some parameters $a$, $\sigma$, $\nu$, $\mu$,
$\lambda$. However, if the intensities imbalance process $r(t)$ is
<<loosened>> and regarded as random, then the final adjustment of
the model reduces to the account of that the parameters $a$,
$\sigma$, $\nu$, $\mu$, $\lambda$ should be regarded as depending on
$r(t)$, so that the final model looks like a a generalized
hyperbolic L{\'e}vy process with random parameters. As this is so,
in accordance with (35), the parameters $a(t)$ and $\sigma(t)$
accumulate the information concerning the current balance between
the sizes and intensities of bid and ask orders, whereas in
accordance with (34) the parameters $\nu(t)$, $\mu(t)$, $\lambda(t)$
are influenced only by the intensities imbalance $r(t)$. The
prediction of the behavior of statistical regularities of this
process and the corresponding risks reduces to the analysis of the
trajectory of a point in a five-dimensional parametric space. For
this purpose one can use multivariate autoregressive models.

\renewcommand{\refname}{References}

\end{document}